\magnification=1200

\def\ni{\noindent}
\def\v{\vskip .3cm}
\def\la{\lambda}
\def\go{{\goth g}}

\def\N{{\bf N}}
\def\Z{{\bf Z}}

\def\sa{\sigma_1}
\def\sb{\sigma_2}
\def\d{\delta}
\def\ds{\delta^{1/2}}
\def\g{\gamma}
\def\span{{\rm span}}
\def\Ca{{\cal C}}

\def\1{1\!\!1}
\def\End{{\rm End}}
\def\Hom{{\rm Hom}}
\def\vl{\vskip .1cm}

\catcode`\@=11    
\def\relaxnext@{\let\next\relax}

\ifx\eulerloaded@AmS\relax\catcode`\@=\active   
 \endinput\else\let\eulerloaded@AmS\relax\fi    
\def\textfont@#1#2{\def#1{\relax\ifmmode
 \Err@{Use \string#1\space only in text}\else#2\fi}}

%
%

\font\teneuf=eufm10

\font\eighteuf=eufm8
\font\seveneuf=eufm7
\font\sixeuf=eufm6  
\font\fiveeuf=eufm5 

%
%

\newfam\euffam      
\textfont\euffam=\teneuf
\scriptfont\euffam=\seveneuf
\scriptscriptfont\euffam=\fiveeuf

%
%

\def\tenpoint{\def\pointsize@{10}%
 \normalbaselineskip12\p@
 \abovedisplayskip12\p@ plus3\p@ minus9\p@
 \belowdisplayskip12\p@ plus3\p@ minus9\p@
 \abovedisplayshortskip\z@ plus3\p@
 \belowdisplayshortskip7\p@ plus3\p@ minus4\p@
 \textfont@\rm\tenrm
 \textfont@\it\tenit
 \textfont@\sl\tensl
 \textfont@\bf\tenbf
 \textfont@\smc\tensmc
 \textfont@\eu\teneuf
 \ifsyntax@\def\big##1{{\hbox{$\left##1\right.$}}}\else
 \let\big\tenbig@
 \textfont\z@\tenrm  \scriptfont\z@\sevenrm  \scriptscriptfont\z@\fiverm
 \textfont\@ne\teni  \scriptfont\@ne\seveni  \scriptscriptfont\@ne\fivei
 \textfont\tw@\tensy \scriptfont\tw@\sevensy \scriptscriptfont\tw@\fivesy
 \textfont\thr@@\tenex \scriptfont\thr@@\tenex \scriptscriptfont\thr@@\tenex
 \textfont\itfam\tenit
 \textfont\slfam\tensl
 \textfont\bffam\tenbf \scriptfont\bffam\sevenbf
  \scriptscriptfont\bffam\fivebf
 \textfont\euffam=\teneuf            
  \scriptfont\euffam=\seveneuf       
  \scriptscriptfont\euffam=\fiveeuf  
 \fi
 \setbox\strutbox\hbox{\vrule height8.5\p@ depth3.5\p@ width\z@}%
 \setbox\strutbox@\hbox{\vrule height8\p@ depth3\p@ width\z@}%
 \normalbaselines\tenrm\ex@=.2326ex}

%
%

\def\eightpoint{\def\pointsize@{8}%
 \normalbaselineskip10\p@
 \abovedisplayskip10\p@ plus2.4\p@ minus7.2\p@
 \belowdisplayskip10\p@ plus2.4\p@ minus7.2\p@
 \abovedisplayshortskip\z@ plus2.4\p@
 \belowdisplayshortskip5.6\p@ plus2.4\p@ minus3.2\p@
 \textfont@\rm\eightrm
 \textfont@\it\eightit
 \textfont@\sl\eightsl
 \textfont@\bf\eightbf
 \textfont@\eu\eighteuf
 \ifsyntax@\def\big##1{{\hbox{$\left##1\right.$}}}\else
 \let\big\eightbig@
 \textfont\z@\eightrm \scriptfont\z@\sixrm  \scriptscriptfont\z@\fiverm
 \textfont\@ne\eighti \scriptfont\@ne\sixi  \scriptscriptfont\@ne\fivei
 \textfont\tw@\eightsy \scriptfont\tw@\sixsy \scriptscriptfont\tw@\fivesy
 \textfont\thr@@\tenex \scriptfont\thr@@\tenex \scriptscriptfont\thr@@\tenex
 \textfont\itfam\eightit
 \textfont\slfam\eightsl
 \textfont\bffam\eightbf \scriptfont\bffam\sixbf
   \scriptscriptfont\bffam\fivebf
 \textfont\euffam=\eighteuf           
  \scriptfont\euffam=\sixeuf          
  \scriptscriptfont\euffam=\fiveeuf   
 \fi
 \setbox\strutbox\hbox{\vrule height7\p@ depth3\p@ width\z@}%
 \setbox\strutbox@\hbox{\vrule height6.5\p@ depth2.5\p@ width\z@}%
 \normalbaselines\eightrm\ex@=.2326ex}

%
%

\def\goth{\relaxnext@\ifmmode\let\next\goth@\else
 \def\next{\Err@{Use \string\goth\space only in math mode}}\fi\next}
\def\goth@#1{{\goth@@{#1}}}
\def\goth@@#1{\fam\euffam#1}

\textfont@\eu\teneuf

\centerline{REPRESENTATIONS OF THE BRAID GROUP B$_3$ AND OF SL(2,{\bf Z})}
\v
\centerline{Imre Tuba and Hans Wenzl
\footnote{$^1$} {Supported in part by NSF grant \# DMS 97-06839}}
\v
\v
We give a complete classification of all simple representations of
$B_3$ for dimension $d\leq 5$ over an algebraically closed field $K$
of any characteristic.
To describe our result in detail, recall that $B_3$ is given by
generators $\sigma_1$ and $\sigma_2$ which satisfy the relation
$\sa\sb\sa=\sb\sa\sb$. Moreover, it is well-known that the
center of $B_3$ is generated by $\zeta=(\sa\sb)^3$.
It is easy to see that $\zeta$ acts on a simple $d$-dimensional
$B_3$- module via the scalar $\delta$ which satisfies the equation
$\delta^d=\det(A)^{6}$, where $A$ is the linear endomorphism
via which $\sa$ acts on $V$.
Our main result states that the eigenvalues of $A$ and the scalar
$\delta$ completely determine a  simple representation of $B_3$
for dimension $\leq 5$, up to equivalence; for $d\leq 3$
it is uniquely determined by the eigenvalues of $A$.
Moreover, such simple representations exist
if and only if the eigenvalues do not belong to the zero set of
certain polynomials in the eigenvalues and $\delta$ which we list
explicitly (see Proposition 2.8, Section 2.10 and 2.11, Remark 4). 


One of the motivations for this paper was studying braided 
tensor categories
by analysing  braid representations. The categories under consideration
have a Grothendieck semiring isomorphic to the one of a semisimple
Lie group. It turns out that in this context
it suffices to classify representations
of $B_3$ up to dimension 5. Indeed, as an application,
we obtain uniform formulas for
the categorical dimensions of objects in the second tensor
power of the adjoint representation in
braided tensor categories closely related to the conjectured
series of exceptional Lie algebras, as proposed by Deligne and Vogel.

Our approach in this paper is quite elementary: we first show that
assuming a certain triangular form of the matrices $A$ and $B$
of the generators
of $B_3$, the braid relation reduces to checking the values
of certain coefficients of the matrix $BA$. We then show that  for dimension
$d\leq 5$ one can always assume such a triangular form, and that
the matrix coefficients are determined by the eigenvalues and $\delta$,
up to certain renormalizations. These explicit representations are then
used to compute polynomials which determine whether the representation is
simple or not. Moreover, we show in the last section how these polynomials
can be used to compute categorical dimensions of objects in braided
categories, and we give explicit formulas for them.

After having completed this research we learned of other approaches
towards classifying braid representations, 
notably via local systems by N. Katz,  and via quivers.
These approaches as well as the case with $d\geq 6$ are discussed in 
the end of chapter 2.

We would like to thank M. Hunziker, N. Wallach, A. Wassermann,
B. Westbury and, in particular, P. Deligne for useful discussions.

\v\ni
{\bf 1.1.} Let $B_3$ be Artin's braid group, given by generators $\sa$ and 
$\sb$ and the relation $\sa\sb\sa=\sb\sa\sb$. It is well-known that
$B_3$ maps surjectively onto $SL(2,\Z)$, e.g. via the map
$$\sa\ \mapsto 
\ \left(\matrix{1&1\cr 0&1\cr}\right),
\quad \quad 
\sb\ \mapsto 
\ \left(\matrix{1&0\cr -1&1\cr}\right).$$
It is easy to check that this homomorphism maps
$$\sa\sb\sa=\sb\sa\sb\ \mapsto\ S\ =\ \left(\matrix{0&1\cr -1&0\cr}\right)
\quad {\rm and}\quad
\sa\sb\ \mapsto 
\ U\ =\left(\matrix{1&1\cr -1&0\cr}\right).$$
Moreover, the center of $SL(2,\Z)$ is equal to $\pm 1$, and the corresponding
elements $\bar S$ and $\bar U$ in the quotient $PSL(2,\Z)$ have orders
2 and 3 respectively. It is known that $PSL(2,\Z)$ is isomorphic to the
free product $\Z_2 *\Z_3$ of a cyclic group of order 2 with a cyclic group
of order 3, and that the isomorphism can be chosen such that $\bar S$ and
$\bar U$ are the generating elements of these cyclic groups. Finally,
the center $Z$ of $B_3$ is generated by $(\sa\sb)^3$, and 
$B_3/Z\cong PSL(2,\Z)$.
\v\ni
{\bf 1.2.} In the following we assume $A$ and $B$ to be invertible $d\times d$
matrices over the algebraically closed field
$K$ satisfying $ABA=BAB$. Hence $\sa\mapsto A$, $\sb\mapsto B$ defines
a representation of $B_3$. Let $\la_1,\ ...\ \la_d$ be
the eigenvalues of $A$ (with $\la_i$ not necessarily distinct from
$\la_j$ for $i\neq j$) and define $\delta=det(A)^{6/d}=
(\prod_{i=1}^d\la_i)^{6/d}$.
\v\ni
{\bf Lemma} (a) Conjugation via $ABA$ maps $A$ to $B$ and $B$ to $A$.

(b) $ABA(AB)^{-1}=B$ and $BAB(BA)^{-1}+A$

(c) $(ABA)^2=\delta 1$, and $(ABA)^{-1}=\delta^{-1}ABA$ if $A$ and $B$
generate the full $d\times d$ matrix ring.

(d) For any choice $\delta^{1/6}$ of a sixth root of $\delta$ the map
$\sa\mapsto \delta^{-1/6}A$, $\sb\mapsto \delta^{-1/6}B$ defines
a representation of $B_3$ and of $PSL(2,\Z)\cong B_3/Z$.

(e) If $\{ e_1,e_2,\ ...\ e_d\}$ is a basis of eigenvectors of $A$,
and $b_i=ABAe_i$, then 
$\{ b_1,b_2,\ ...\ b_d\}$ is a basis of eigenvectors of $B$.
\v
$Proof$ It follows from the braid relation that $BAB^{-1}=A^{-1}BA$,
from which one easily deduces claims (a) and (b). For (c), observe that 
$(\sa\sb\sa)^2=(\sa\sb)^3$ is in the center of $B_3$. As we assume
an irreducible representation of $B_3$, the center
has to act via a scalar matrix,
say $\la 1$. Comparing determinants (using $det(A)=det(B)$, 
we get $\la^d=\det(A)^6$. Statement (e) follows from (a); (d) is clear.

\v\ni
{\bf 1.3.} We say that the $d\times d$ matrices
$A$ and $B$ are in {\it ordered triangular} form if
$A$ is an upper triangular matrix with eigenvalue $\la_i$ as
$i$-th diagonal entry, and $B$ is a lower triangular matrix with
eigenvalue $\la_{d+1-i}$ as $i$-th diagonal entry.
\v\ni
{\bf Lemma} Assume $A$ and $B$ are in ordered triangular form, satisfying
$ABA=BAB$. Then

(a) Any $d\times d$ matrix $F=(f_{ij})$ such that $FBF^{-1}=A$ 
has matrix entries $f_{ij}=0$
for $i+j>d+1$. This holds, in particular, for $F=BA$.

(b)  Any matrix $d\times d$ $E=(e_{ij})$ such that $EAE^{-1}=B$ 
has matrix entries $e_{ij}=0$
for $i+j<d+1$. This holds, in particular, for $E=AB$.

\v
$Proof.$ Let $\{ w_1,w_2,\ ...\ w_d\}$ be the standard basis for $K^d$.
By assumption, the eigenspace of $A$ corresponding to eigenvalues
$\la_1,\ ...\ \la_r$ is equal to span$\{ w_1,w_2,\ ...\ w_r\}$, for
$r=1,2,\ ...\ d$, while the corresponding eigenspace for $B$ is
equal to span $\{w_{\bar 1},w_{\bar 2},\  ...\ w_{\bar r}\}$,
where $\bar i=d+1-i$. Hence $F$ maps
span$\{ w_1,w_2,\ ...\ w_r\}$ onto
span$\{ w_{\bar 1},w_{\bar 2}\ ...\ w_{\bar r}\}$,
from which we deduce the triangular shape
of $F$. The second statement in (a) follows from Lemma 1.2(b).
Statement (b) is proved similarly.
\v\ni
{\bf Proposition} Assume $A$ and $B$ are in ordered triangular form, 
satisfying $ABA=BAB$. Let $\bar i=d+1-i$.
After rescaling our basis vectors, if necessary, we can assume
\vskip .1cm
(a)\quad   $ABA$ is a skew diagonal matrix, with $(ABA)_{i,\bar i}=
(-1)^{i+1}\delta^{1/2}$,
\vskip .1cm
(b)\quad  $b_{ij}=(-1)^{i+j}a_{\bar i,\bar j}$,
\vskip .1cm
(c)\quad  $\sum_{k=1}^i(-1)^{k+i}a_{\bar i,\bar k}a_{k,\bar i}=
(BA)_{i,\bar i}=(-1)^{i+1}\delta^{1/2}/\la_i$, 
\vskip .1cm
(d)\quad  If $i+j>d+1$, $\sum_{k=1}^{min(i,j)}(-1)^{k+i}
a_{\bar i,\bar  k}a_{k,j}=(BA)_{ij}=0$

(e)\quad Statements (a)-(d) also hold for the coefficients of the
matrices $A'=DAD^{-1}$ and $B'=DBD^{-1}$ if $D$ is a diagonal matrix
with $d_{ii}=d_{\bar i\bar i}\neq 0$. In particular, if the entries
$a_{id}$ for $1<i\leq [(d+1)/2]$ are nonzero, we can choose arbitrary
nonzero numbers for them.
\v
$Proof.$ Statements (a) and (b) follow from the previous Lemma and 
Lemma 1.2(a) and (c).
Statements (c) and (d) follow from explicit matrix computations,
exploiting  the triangular shapes of $A$, $B$
and $BA$, and using statement (a). For (e), observe that conjugation
by $D$ does not change $ABA$, and also leaves $A'$ and $B'$ in ordered
triangular form.
One can now prove statements (a)-(d) for $A'$ and $B'$ as 
we did before for $A$ and $B$.
\v\ni
{\bf Corollary} (a) $a_{1d}=\delta^{1/2}/\la_1\la_d$,,

(b) $a_{2,d-1}=(\la_2a_{d-1,d}a_{1,d-1}-\delta^{1/2})/\la_2\la_{d-1}$,
\v
$Proof.$ By Proposition, (c), we have $\delta^{1/2}/\la_1=(BA)_{1d}$,
with the latter being equal to $a_{1d}\la_d$. Statement (b) follows
similarly, using Proposition (c). 
\v\ni
{\bf 1.4. Proposition } Let $A$ be an upper triangular matrix with
eigenvalues $\la_1, \la_2,\ ...\ \la_d$ down the diagonal, and let
$S$ be a skew-diagonal matrix with $S^2=c1$ for some constant $c$
such that $B=SAS^{-1}$.
Moreover, assume

(a) $(BA)_{ij}=0$ for $i+j>d+1$,

(b) $\la_i(BA)_{i\bar i}=cs_{i,\bar i}$

\ni Then $A$ and $B$ satisfy the braid relations.

\v
$Proof.$ Observe that $S^{-1}$ is a scalar multiple of $S$. 
Hence $B=S^{-1}AS=SAS^{-1}$.
We get from this $AB=S(BA)S^{-1}$, hence $(AB)_{ij}=0$ for
$i+j<d+1$. Exploiting the triangular shapes of the matrices, we get
that $A(BA)$ is upper skew-triangular and $(AB)A$ is lower skew-triangular.
Hence $ABA$ is of skew-diagonal shape, with 
$(ABA)_{i,\bar i}=\la_i (BA)_{i,\bar i}=cs_{i\bar i}$,  
and $ABA$ commutes with $S$. Hence 
$$BAB=S(ABA)S^{-1}=ABA.$$
\v\ni
{\bf 1.5. An example} Let $V$ be a $(d+1)$-dimensional vector space,
with a basis labeled by $0,1,\ ...\ d$, and let $\la_0,\la_1,\ ...\ \la_d$
be parameters satisfying
$\la_i\la_{d-i}=c$ for a fixed constant $c$. Only for this
subsection we define $\bar i=d-i$. Then we get a $[{d+1\over 2}]$
parameter family of representations of $B_3$ via the matrices
$$A=\left(\left(\matrix{\bar i\cr \bar j\cr}\right)\la_j\right)_{ij},
\quad
B=\left((-1)^{i+j}\left(\matrix{i\cr j\cr}\right)\la_{\bar i}\right)_{ij}.$$
To prove this, it satisfies to  check the conditions of Proposition 5,
with $S$ being the skew-diagonal matrix defined by
$s_{ij}=(\delta_{i,\bar j}(-1)^i\la_{\bar i})$. This can be fairly easily
checked, using the identity
$$\sum_{k=0}^d(-1)^{i+k}\left(\matrix{i\cr k\cr}\right)
\left(\matrix{\bar k\cr \bar j\cr}\right)=
(-1)^i\left(\matrix{d-i\cr d-j\cr}\right)=
(-1)^i\left(\matrix{\bar i\cr \bar j\cr}\right).$$
This identity is well-known, and can be easily proved by induction on $i$.
\v\ni
{\bf 2.} Representations of $B_3$ of dimension $d\leq 5$.
\v\ni
{\bf 2.1}
In the following let $V$ be a $B_3$-module over the field $K$. Let
$\{ e_1,e_2,\ ...\ e_d\}$  be a basis of generalized eigenvectors of $A$, and
let $\{ b_1,b_2,\ ...\ b_d\}$ be a basis of generalized eigenvectors of $B$
with $b_i=ABAe_i$ for $i=1,2,\ ...\ d$. We shall always assume
that if $A$ and $B$ have Jordan blocks,
the labeling is chosen such that they are upper triangular. 
As a consequence of 
our assumptions and Lemma 1.2, there exist scalars $c_{ij}$,
$1\leq i,j\leq d$ such that
$$b_i=\sum_{j=1}^d c_{ji}e_j
\quad {\rm and}\quad e_i=\delta^{-1}\sum_{j=1}^d c_{ji}b_j.$$
As $(ABA)^2=\delta 1$, the matrix $ABA$ is diagonalizable. We will usually
refer to an eigenvalue of $ABA$ just as $\delta^{1/2}$, not specifying
which root we choose.
We  shall also consider the subspaces
$$W=\span\{ e_1,e_2,\ ...\ e_{d-1}\}\ \cap\ \span\{ b_1,b_2,\ ...\ b_{d-1}\}$$
$$ W'=\bigcap_{i\in\N}
(AB)^i(\span\{ e_1,e_2,\ ...\ e_{d-1}\}).$$
Observe that $W'\subset W$, as $ABe_i=ABAA^{-1}e_i=\la_i^{-1}b_i$. Moreover,
$W$ is invariant under $ABA$ and $W'$ is invariant under $AB$.
Also observe that $W$ has codimension $\leq 2$ and $W'$ has codimension
$\leq 3$, being the intersection of 2 resp. 3 subspaces of $V$
with codimension 1. We shall need the following simple observations:

(i) Each of the sets $\{ \sa, \sb\}$,   $\{ \sa, \sa\sb\}$ or 
$\{ \sa\sb\sa, \sa\sb\}$
generates $B_3$.

(ii) Let $I$ be a subset of $\{ 1,2,\ ...,\ d\}$. If $V'=\span\{ e_i,\ i\in I\}
=\span\{ b_i,\ i\in I\}$, then $V'$ is a $B_3$-submodule.

(iii) If $W$ has codimension 1, it would coincide with
$\span\{ e_1,e_2,\ ...\ e_{d-1}\}$ and with
$\span\{ b_1,b_2,\ ...\ b_{d-1}\}$, and therefore would be a $B_3$-submodule,
by (ii).

(iv) If both $W$ and $W'$ have codimension 2, they would coincide; this  space
would be invariant under both $ABA$ and $AB$, from which one easily 
deduces that it is a $B_3$-submodule, by (i).
\v\ni
{\bf 2.2. Proposition} Let $V$ be a simple $B_3$-module of dimension
$d\leq 5$. Then an eigenvector of $A$, say $e_i$
can not be contained in a proper subspace of
$V$ which is invariant  under $B$ and contains $b_i=ABAe_i$.
\v
$Proof.$ We choose the labeling of (generalized) eigenvectors so that
$e_1$ is contained in a $B$-invariant subspace spanned by generalized
eigenvectors $b_1$, ..., $b_{d-1}$.
Let $W$ and $W'$ be as in 2.1.
By 2.1(iii) and (iv), the claim follows immediately unless
$W$ has codimension 2 and $W'$ has codimension 3.
As $W$ has at least dimension 1, containing $e_1$, we get the claim for $d=2$
immediately.

$d=3$: Here $W$ has to be 1-dimensional, containing
both $e_1$ and $b_1$, i.e. it is an eigenspace of both $A$ and $B$.

$d=4$: By assumption, $W$ contains both $e_1$ and $b_1$ and has dimension
2. If $e_1$ and $b_1$ are  linearly dependent, the claim follows from
(2.1)(ii). Hence we can assume them to be a basis for $W$.
As dim $W'=1$, there exists an eigenvector $v=\alpha e_1+\beta b_1$
of $AB$ in $W'$. We can assume $\beta\neq 0$, as otherwise $v$ would
be an eigenvector of both $A$ and $AB$, from which the claim
would follow by 2.1(i). But then
$$\mu(\alpha e_1+\beta b_1)=AB(\alpha e_1+\beta b_1)=
\la_1^{-1}b_1+\beta\la_1Ab_1,$$
where $\mu$ is the eigenvalue of $v$. We deduce from this $Ab_1\in W$,
and also $\delta Be_1=BABAb_1\in BAB(W)=W$.
Hence $W$ is invariant under both $A$ and $B$.

$d=5:$ We can assume $W$ to have a basis $\{ e_1, b_1, x\}$, with $x$
an eigenvector of $ABA$ with eigenvalue $\ds$ (see 2.1). 
As $AB$ is diagonalizable,
we can choose a basis of eigenvectors 
$v_i=\alpha^{(i)}e_1+\beta^{(i)}b_1
+\gamma^{(i)}x$, $i=1,2$ for $W'$.
We have $ABv_i=w_i+\alpha^{(i)}\la_1^{-1}b_1$, with 
$w_i=\la_1\beta^{(i)}Ab_1+\gamma^{(i)}ABx$, for $i=1,2$. As $ABv_i$ and 
$b_1$ are in $W$, so is $w_i$, for $i=1,2$.

$Case\ 1:$ If $w_1$ and $w_2$ are linearly dependent, then there
exist scalars $\nu_1$ and $\nu_2$ such that $\nu_1v_1+\nu_2v_2\neq 0$
and $\nu_1w_1+\nu_2w_2
=0$. But then $AB(\nu_1v_1+\nu_2v_2)\in W'$ is a nonzero multiple of $b_1$;
in particular, $b_1\in W'$. But then also
$$\delta\la_1e_1= \delta Ae_1=A(BAB)b_1=(AB)^2b_1\in W'.$$
Hence $W'$ is spanned by $e_1$ and $b_1$, and therefore is also invariant
under $ABA$.

$Case\ 2:$ If $w_1$ and $w_2$ are linearly independent, then also
$ABx$ and $Ab_1$ are  in $W$. We conclude $Be_1\in W$ as in the case $d=4$.
Moreover, also $\ds B^{-1}x=B^{-1}BABx=ABx\in W$, i.e.
$B^{-1}x$ is a linear combination of $e_1,b_1$ and $x$. If the coefficient
of $x$ is not equal to 0 in this linear combination, we can multiply it
by $B$ and solve for $Bx$, which shows that it is in $W$. In this
case, $W$ is invariant under both $B$ and $BAB$.
If $B^{-1}x$ is a linear combination of only $e_1$ and $b_1$, we can multiply
this linear combination by $AB$ to obtain for $Ax$ a linear combination
in $b_1$ and $Ab_1\in W$. Hence $W$ is invariant under $A$ and $ABA$.
The claim follows, using 2.1(i).
\v\ni
{\bf Corollary} Let $V$ be a simple $B_3$-module with dimension $d\leq 5$.
Then the minimal polynomial of $A$ coincides with its characteristic 
polynomial.
\v
$Proof.$ Let $b\in V$, let $\tilde S = \span\{ A^ib,\ i=0,1,\ ...\ \}$,
and let $\la$ be an eigenvalue of $A$. Then it is well-known
that the intersection of the eigenspace of the eigenvalue $\la$ (for $A$)
with $\tilde S$ has at most dimension 1.
To see this directly, 
let $p_\la$ be the projection onto the generalized eigenspace of $A$
for the eigenvalue $\la$, with the kernel being the direct sum of the
generalized eigenspaces for the other eigenvalues. Then we get
$$p_\la\tilde S=\span\{ A^ip_\la b, i=0,1,\ ...\ \} = 
\span\{ (A-\la)^ip_\la b, i=0,1,\ ...\ \}.$$
The claim now follows from the fact that the nonzero elements in
$\{ (A-\la)^ip_\la b, i=0,1,\ ...\ \}$ are linearly independent,
and $A$ acts as a full Jordan block on their span.

If the minimal polynomial of $A$ does not coincide with the
characteristic polynomial, there exists an eigenvalue, say $\la_1$,
whose eigenvectors span a subspace $E$ of dimension at least 2.
If there is an eigenvalue distinct from $\la_1$,  say $\la_j$, 
pick an eigenvector
$b_j$ of $B$ belonging to $\la_j$. Then 
$\span\{ A^sb_j, s=0,1,\ ...\ d-2\}$ forms an $A$-invariant subspace
$S'$ such that $\dim S'\cap E\leq 1$, as proved in the last
paragraph. Let $S$ be the subspace generated by
$S'$ and the eigenvectors of $A$ with eigenvalue $\la_j$. Then
also $\dim S\cap E\leq 1$, i.e. $S$ is a proper $A$-invariant subspace of
$V$ which contains both $b_j$ and $BAb_j$. Hence $V$ could not be
a simple $B_3$-module.

If $A$ only has one eigenvalue, $\dim E\geq 2$ implies that we have
at least 2 different Jordan blocks. Let $b_j$ be an eigenvector of $B$
belonging to a block of maximum length. Let $S'=
\span\{ A^sb_j, s=0,1,\ ...\ d-2\}$, which is a proper subspace
of $V$ by assumption. If $S'$ did not contain $e_j=BAb_j$,
it would not contain the whole Jordan block of $e_j$ either.
Hence the space $S$ spanned by $S'$ and $\{ e_j\}$ is a proper $A$-invariant
subspace containing both $e_j$ and $b_j$, and hence is a $B_3$-submodule.
If $S'$ does contain $e_j$, it would be a proper $B_3$-submodule by
the same argument.
\v\ni
{\bf 2.3 Lemma} Let $e_1$ and $e_2$ be either 2 eigenvectors of $A$ with
 eigenvalues $\la_1$ and $\la_2$, or let them be generalized
eigenvectors with $Ae_1=\la e_1$ and $Ae_2=\la e_2+e_1$.

(a) Let $V'=\span\{ e_1,e_2,b_1,b_2\}$. If $\dim V'\leq 3$, then there exists
a nonzero $B_3$-invariant subspace of $V$ of dimension $\leq 3$.

(b) Let $d=5$ and let $V''=\span\{ e_1, e_2, e_3, b_1, b_2\}$, with
$e_3$ a generalized eigenvector belonging to $e_1$ or $e_2$, or
an eigenvector.
If $\dim V''\leq 4$, it contains a nonzero
$B_3$-invariant subspace.
\v
$Proof.$ The claim in (a) follows from 2.1(ii) if dim $V'=2$. Hence
we can assume  dim $V'=3$.
Let $W=\span\{ e_1,e_2\} \cap \span\{ b_1,b_2\}$. Then 
 $\dim\ W=1$, as otherwise $\{ e_1, e_2, b_1, b_2\}$ would be  linearly
independent. In particular, $W$ is spanned by an
eigenvector $w=\alpha_1e_1+\alpha_2e_2$ of $ABA$. 
Let us assume first that both $e_1$ and $e_2$ are eigenvectors.
If $\la_1=\la_2$,
$w$ would be a common eigenvector of both $A$ and $ABA$, and we are done.
Hence we can assume $\la_1\neq \la_2$.
We can also assume $\alpha_2\neq 0$ in the expression for $w$, as
otherwise $w$ would be an eigenvector of both $A$ and $ABA$.
Moreover, we also have $\ds w=ABAw=\alpha_1b_1+\alpha_2b_2$. We now
compute
$$ABw=B^{-1}BAB(\alpha_1e_1+\alpha_2e_2)
=\alpha_1\delta^{1/2}w+\alpha_2(\la_2^{-1}-\la_1^{-1})b_2,$$
$$(AB)^2w=\delta^{1/2}A(\alpha_1e_1+\alpha_2e_2)=
\delta^{1/2}(\la_1w+\alpha_2(\la_2-\la_1)e_2).$$
The coefficients
of $e_2$ and $b_2$ in the expressions above are nonzero. Hence the vectors
$\{e_2,b_2,w\}$  and
$\{ w,\ ABw,\ (AB)^2w\}$
are bases for the same subspace.
This proves that it
is both invariant under $AB$ (second basis, as $(AB)^3w=\d w$)
 and under $ABA$ (first basis). This proves statement (a) if both
$e_1$ and $e_2$ are eigenvectors. If $e_2$ is a generalized eigenvector,
we compute
$$ABw=\ds \la^{-1}w-\ds\alpha_2\la^{-2}b_1,\quad (AB)^2w=
\ds \la w+\ds \alpha_2e_1.$$
One shows as before that $\span\{ w,\ ABw,\ (AB)^2w\}=\span\{e_1,b_1,w\}$, 
from which one deduces the claim.

Observe that $V'\subset V''$. Statement (b) follows immediately from (a)
if $\dim V'<4$. Hence we can assume $\dim V'=4$ and $V'=V''$, which therefore
is invariant under $ABA$. Moreover, the space $U=\span\{ e_1,e_2,e_3\}\cap
\span\{ b_1,b_2,b_3\}$ has at least dimension 2, being the intersection
of 2 subspaces of $V'$ with codimension 1; we can assume
$\dim U=2$ as otherwise $U$  would be an invariant
subspace. Hence it contains 
2 linearly independent eigenvectors $w_1$ and $w_2$ of $ABA$. Moreover,
if $w_1=\sum \alpha_je_j$, we can assume $\alpha_3\neq 0$; otherwise
$w_1\in \span\{ e_1,e_2\}\cap \span\{ b_1,b_2\}$ which is 0 (the vectors
$  e_1,e_2,  b_1,b_2$ are a basis for $V'$).
Now the equality $ABAw_1=\ds w_1$ implies 
$$BAw=B(\sum_{j=1}^3\alpha_j\la_je_j)=\sum_{j=1}^3\ds\alpha_j\la_j^{-1}e_j
= \ds A^{-1}w.\eqno(*)$$
We want to show that we can assume $e_3\not\in U$. If $e_3$ is a generalized
eigenvector belonging to $e_1$ with $e_3\in U$,
we can replace it by $e_3+\tau e_1$ for any scalar $\tau\in K$. This new vector
would only be in $U$ for any choice of $\tau$ if $e_1\in U$. But then
$U$ would be spanned by $e_1$ and $e_3$, and would be a $B_3$-submodule,
being invariant under $A$ and $ABA$. 
If $e_3$ is an eigenvector in $U\subset \span\{ b_1,b_2, b_3\}$,
the claim follows from Prop. 2.2.

So we can assume
that $X=\span \{e_3, w_1, w_2\}$ coincides with
$\span\{ e_1,e_2,e_3\}$. 
If $e_4$ and $e_5$ were eigenvectors with same eigenvalue $\la$,
then either $ABAe_3\in X$ (which would make $X$  invariant under
both $A$ and $ABA$), or $\span(\{ ABAe_3\}\cup X)$ would be equal to
$\span(\{ \tilde e_4\}\cup X)$, where
$\tilde e_4$ is a  linear combination of  $e_4$ and $e_5$, and, 
in particular, it is an eigenvector of $A$. Using these 2 different
spanning sets, we see that we would obtain a 4-dimensional subspace
invariant under both $A$ and $ABA$.

Hence we can assume that $V/X$ has  at most 2 $A$-invariant 1-dimensional
subspaces, and, similarly $V/\span\{ b_1, b_2, b_3\}$ has at most
2 $B$-invariant 1-dimensional  subspaces.
Expanding the vectors in the  equation $(*)$
as a linear combination of $ \{e_3, w_1, w_2\}$, 
and  observing that $w_1,w_2\in 
\span\{ b_1,b_2,b_3\}$, we obtain that $Be_3$ is congruent to a multiple
of $e_3$ modulo $\span\{ b_1,b_2,b_3\}$. But this would imply that
$e_3$ together with $ b_1,b_2,b_3$  spans a $B$-invariant subspace of
dimension $\leq 4$. Statement (b) now follows from Proposition 2.2,
if $e_3$ is an eigenvector. If $e_3$ is a generalized eigenvector
belonging to, say $e_1$, we could show as before that also 
$e_3'=e_3+\tau e_1$ is in a 4-dimensional $B$-invariant subspace
containing $b_1, b_2$ and $b_3$. As this is true for any $\tau\in K$,
and there are  only at most 2 such subspaces, we obtain that also
$e_1$ itself is in a $B$-invariant subspace containing $b_1$.
The claim follows from Proposition 2.2.
\v\ni
{\bf 2.4. Proposition} Let $V$ be a simple $B_3$-module with dimension
$\leq 5$. Then there exists a basis of $V$ with respect to which
$A$ and $B$ act in ordered triangular form (see 1.3). 
Moreover, this is possible for any labeling of the generalized eigenvectors
as long as $A$ appears in upper triangular Jordan form.

\v
$Proof.$
To construct such a basis,
we can always assume at least one eigenvector $e_1$ for $A$
and $b_1$ for $B$, and at least one more (generalized) eigenvector
$e_2$ as in Lemma 2.3. We define
$w_1=e_1$ and $w_d=b_1$. 
For $d\geq 3$, observe that $c_{d1}\neq 0$ by Proposition 2.2,
with $c_{ij}$ as in 1.3, for $1\leq i,j\leq d$.
Hence we can define $w_2=e_2-(c_{d2}/c_{d1})e_1$, which is both in
$\span\{ e_1, e_2\}$ and in $\span\{ b_1,\ ...,\ b_{d-1}\}$.
For $d\geq 4$, define $w_{d-1}=ABAw_2$, which is in
$\span\{ b_1, b_2\}$ and in $\span\{ e_1,\ ...,\ e_{d-1}\}$.
Finally, if $d=5$, we can express $b_3$ as a linear combination
of $\{ e_1, e_2, e_3, b_1, b_2\}$, by Lemma 2.3. Hence there exist scalars
 $\alpha_1$ and $\alpha_2$ such that
$w_3=b_3-\alpha_1b_1-\alpha_2b_2$ is in
$\span\{ e_1, e_2, e_3\}\cap\span\{ b_1, b_2, b_3\}$. 
By construction, it follows that
$$w_i\in\span\{ e_1, e_2,\ ...\  e_i\}\cap\span\{ b_{1},b_{2},
\ ...\  b_{d+1-i}\}\quad {\rm for\ }1\leq i\leq 3\eqno(*)$$
It follows from $(*)$ 
that $\span\{w_1,w_2,\ ...\ w_i\}\subset\span\{e_1,e_2,\ ...\ e_i\}$
and  that 

\ni
$\span\{w_{d+1-i},w_{d+2-i},\ ...\ w_d\}
\subset\span\{b_1,b_2,\ ...\ b_i\}$. Let us check equality for
$d=5$. This follows by construction for $i\leq 3$. 
But then $\span\{ w_1,w_2,w_3\}=\span\{ e_1,e_2,e_3\}$ and
$\span\{w_4,w_5\}=\span\{ b_1,b_2\}$. The linear independence of
the $w_i$'s now follows from Lemma 2.3(b). Hence the inclusions
below $(*)$ actually are equalities, from which easily
deduces the triangular forms
of $A$ and $B$ by induction on $i$. The cases $d<5$ are similar and
easier to check.

\v\ni
{\bf 2.5. Proposition} Let $V$ be a simple $B_3$-module of dimension $d$
with $d=2,3$. Then there exists a basis for $V$ for which 
$A$ and $B$ acts via the matrices
$$A=\left(\matrix{\la_1&\la_1\cr 0&\la_2\cr}\right),
\quad\quad
B=\left(\matrix{\la_2&0\cr -\la_2&\la_1\cr}\right)
\quad {\rm for\ } d=2$$
\v
$$A=\left(\matrix{\la_1&\la_1\la_3\la_2^{-1}+\la_2&\la_2\cr
            0&\la_2&\la_2\cr 0&0&\la_3\cr}\right),
\quad\quad
B=\left(\matrix{\la_3&0&0\cr
        -\la_2&\la_2&0\cr \la_2&-\la_1\la_3\la_2^{-1}-\la_2&\la_1\cr}\right)
\quad {\rm for\ } d=3. $$
\v
$Proof.$ For $d=2$: By Proposition 2.4, we can assume $A$ and $B$ in ordered
triangular form, with only the nonzero off-diagonal entries to be
computed.  Rescaling one of the basis vectors, we can  assume
$a_{12}=\la_1$. We obtain $b_{21}=-\la_2$ from the braid relation
$ABA=BAB$.

For $d=3$: Again, by Proposition 2.4, we only need to compute
the nonzero off-diagonal entries of
$A$ and $B$, which  are in ordered triangular form. By Corollary 1.3,
we can assume $a_{13}=\la_2=b_{31}$. If $a_{23}=0$, then also $b_{21}$,
by Proposition 1.3(b). But then span$\{w_1,w_3\}$ would be a subspace
invariant under both $A$ and $B$, contradicting $V$ being simple.
Hence we can assume $a_{23}=\la_2=-b_{21}$, by Prop. 1.3(e). 
Finally, we compute
from $(BA)_{33}=0$ that $b_{32}=-\la_1\la_3\la_2^{-1}-\la_2$.
\v\ni
{\bf 2.6. Proposition} Let $V$ be a 4-dimensional $B_3$-module, and let
$D =\sqrt{\la_2\la_3/\la_1\la_4}$. Then we can
find a  basis for $V$ with respect to which we get the matrices
$$A=\left(\matrix{\la_1&(1+D^{-1}+D^{-2})\la_2&(1+D^{-1}+D^{-2})\la_3&\la_4\cr
0&\la_2&(1+D^{-1})\la_3&\la_4\cr
0&0&\la_3&\la_4\cr
0&0&0&\la_4\cr}\right)$$
$$B=\left(\matrix{\la_4&0&0&0\cr
-\la_3&\la_3&0&0\cr D\la_2&-(D+1)\la_2&\la_2&0\cr
-D^3\la_1&(D^3+D^2+D)\la_1&-(D^2+D+1)\la_1&\la_1\cr}\right).$$

$Proof.$
For the proof, we first assume $A$ and $B$ in the form of Proposition 1.3.
Then we get from its corollary
$$a_{14}=\d^{1/2}/\la_1\la_4=(\la_1\la_2\la_3\la_4)^{1/4}D,\quad\quad
a_{23}=(a_{34}a_{13}\la_2-\d^{1/2})/\la_2\la_3,$$

$$a_{34}={\la_3\over a_{14}}a_{24}={\la_1\la_3\la_4\over 
(\la_1\la_2\la_3\la_4)^{3/4}}a_{24},
\quad\quad
a_{13}={(\la_1\la_2\la_3\la_4)^{3/4}\over \la_1\la_2\la_4}a_{12},$$
where the last 2 equalities follow from $(BA)_{24}=0$ and $(BA)_{42}=0$.
It follows from the equations above that $a_{34}=0$ if and only if $a_{24}=0$,
and in this case, also $b_{21}=0=b_{31}$, by Prop. 1.3(b).
But then $span\{w_2, w_3\}$ would be a $B_3$-submodule, contradicting
simplicity of $V$. 
By Proposition 1.3,(e) we can choose 
$a_{24}=\la_4$, from which one deduces
$$ a_{34}= {\la_1\la_3\la_4^2\over
(\la_1\la_2\la_3\la_4)^{3/4}}$$
We get from $(BA)_{34}=0$, using the substitution in Prop. 1.3(b) that
$$a_{23}=(\la_1\la_2\la_3\la_4)^{1/4}(D+1).$$
Similarly, we obtain from $(BA)_{44}=0$, using Prop. 1.3(b) and the results
so far
$$a_{13}=\la_4^{-1}\la_2\la_3(D+1+D^{-1}),
\quad
a_{12}={\la_1\la_2^2\la_3\over(\la_1\la_2\la_3\la_4)^{3/4}}(D+1+D^{-1}).$$
To get $A$ and $B$ into the form as stated, it suffices to conjugate
the matrices via the diagonal matrix diag$(\la_4/a_{i4})_i$.
This shows that there exist at most 2 representations of $B_3$
with prescribed eigenvalues $\la_i$ for $A$, up to conjugation,
depending on the choice of the square root in the expression for $D$.

Observe that we used the equations $(BA)_{ij}=0$ for $j=4$ and $i>1$
for the computation of matrix entries above. It is easy to check
$(BA)_{ij}=0$ for the remaining entries for which $i+j>5$. Hence
condition (a) of Proposition 1.4 is satisfied, and condition (b)
follows from Prop 1.3(d). Hence $A$ and $B$ satisfy the braid relations.

\v\ni
{\bf 2.7. Proposition} Let $V$ be a simple 5-dimensional $B_3$-module.
Then there exists a basis with respect to which $A$ and $B$ act
as ordered triangular matrices.
Moreover, the module is uniquely determined by the eigenvalues of $A$
and the choice $\g$ of a 5th root of det$(A)$, up to equivalence.
\v
$Proof.$ We proceed as in the proof of Proposition 2.5, using
the results of Section 1.3. By Corollary 1.3(a), we get
$$a_{15}={\ds \over \la_1\la_5}={\g^3\over \la_1\la_5}.$$

If $a_{35}=0$, we can conclude from $(BA)_{35}=0$ that 
$a_{34}a_{25}=a_{35}(\la_3+a_{15})=0$. 
If $a_{34}=0$, then also $b_{31}=0=b_{32}$ by Prop. 1.3(b). In this case,
span$\{ w_1,w_2,w_4,w_5\}$ is a $B_3$-submodule.
Hence we can  assume $a_{34}\neq 0$ and $a_{25}=0$. But then
$a_{45}=(\la_4/a_{15})a_{25}=0$ too by Corollary 1.3(b); 
and the matrix entries
$b_{21}, b_{31}$ and $b_{41}$ are also equal to 0, by Prop. 1.3(b). Hence
span$\{ w_1,w_5\}$ is an invariant subspace.

So we can assume $a_{35} \neq 0$.
Assume $a_{45}=0$; then also $a_{25}=(a_{15}/\la_4)a_{45}=0$.
We get from $(BA)_{45}=0$ that
$a_{35}a_{23}=\la_2a_{45}+a_{24}a_{25}-a_{25}a_{15}=0.$ Since $a_{35} \neq 0$,
$a_{23}=0$. It follows that $b_{21},b_{41}$ and $b_{43}$
are all equal to 0, by Prop 1.3(b).
This would entail that span$\{ w_1, w_3, w_5\}$
is a $B_3$-submodule.

Hence we can assume that $a_{i5}\neq 0$ for $i=1,2,\ ...\ 5$.
We can choose $a_{45}=\la_4$ and $a_{35}=a_{15}$, by Prop. 1.3(e).
We now show that the equations in (1.3) completely determine
the other entries of $A$ and $B$. As this does not seem to be
completely straightforward, we include the details for the
interested reader.
Using $(BA)_{25}=0$ and $(BA)_{35}=0$ we get
$$a_{25}=a_{15}={\g^3\over \la_1\la_5}, \quad\quad
a_{34}=a_{15}+\la_3 = {\g^3\over \la_1\la_5} + \la_3.$$
>From $(BA)_{24}=-{\ds \over \la_2} = -{\g^3\over \la_2}$,
$$a_{24}=a_{14}-{\g^3\over \la_2\la_4}.$$
Substituting this into $(BA)_{34}=0$, we find
$$a_{14}=
{a_{34}(\la_2\la_3\la_4+\g^3)\over (a_{34}-a_{35})\la_2\la_4} = 
\left({\la_2\la_4\over \g^2}+1\right)
\left(\la_3+{\g^3\over \la_2\la_4}\right),$$
and hence
$$a_{24}={\g^3\over \la_1\la_5}+\la_3+\g.$$
Using $(BA)_{52}=0$ we obtain
$$a_{12}={\la_2a_{14}\over a_{15}} = \left(1+{\g^2\over
\la_2\la_4}\right)\left(\la_2+{\g^3\over \la_3\la_4} \right).$$
Now $(BA)_{45}=0$ gives us
$$a_{23}=a_{24}-a_{15}+{\la_2\la_4\over a_{15}} = \g + \la_3 + {\g^2\over
\la_3}.$$
Finally, using  $(BA)_{55}=0$, we get
$$a_{13}=\left({\g^2\over \la_3}+\la_3+\g\right)\left(1 + {\la_1\la_5\over 
\g^2}\right).$$
This shows that all entries of $A$ and $B$ are uniquely determined by the
eigenvalues $\la_i$ of $A$ and by a choice of a $5^{th}$ root of
det$(A)$. It is now a straightforward
computation to check the remaining conditions of
Prop. 1.4 to prove that the matrices $A$ and $B$ do indeed define
a representation of $B_3$. Alternatively, the existence question
might be more easily settled using some of the methods discussed in
Section 2.11.

\v\ni
{\bf 2.8} Define for $1\leq i\leq d$ the polynomials $P_r^{(d)}(x)
=\prod_{i\neq r}(x-\la_i)$, where $1\leq i\leq d$ with $i\neq r$.
Observe that these polynomials divide the characteristic polynomial
of $A$. It follows from the corollary of Prop. 2.2 that
$P_r^{(d)}(A)$ is a nonzero rank 1 matrix; if the $d$ eigenvalues of $A$ are
mutually distinct, it is a multiple of the projection onto the eigenspace
of $\la_r$, with kernel being the direct sum of the eigenspaces of the
other eigenvalues.

Now observe that for our braid representations, $A$ and $B$ are matrices
with coefficients in the ring $R$ of Laurent polynomials in  the $\la_i$'s and
$\g$, a $d$-th root of $\det(A)$ for $d=4,5$.
As  also the matrix $P_r^{(d)}(A)P_s^{(d)}(A)P_r^{(d)}(A)$ is
a multiple of  $P_r^{(d)}$, we obtain Laurent
polynomials $Q_{rs}^{(d)}$ in $R$ by
$$P_r^{(d)}(A)P_s^{(d)}(B)P_r^{(d)}(A)=Q_{rs}^{(d)}P_r^{(d)}(A).$$
\v\ni
{\bf Proposition} (a) $P_1^{(d)}(B)P_d^{(d)}(A)=Q_{1d}^{(d)}E_{dd}$,
where $E_{dd}$ is the matrix with a 1 in the  $(dd)$-entry
and zeroes everywhere  else.

(b) If $Q_{ij}^{(d)}=0$, then $P_i^{(d)}(B)P_j^{(d)}(A)=0$.

(c) The matrix $(\prod_{r=i}^d(A-\la_r))$
has nonzero entries in its last column at most
in the 1-st until $i-1$-st row.

(d) The Laurent polynomials $Q_{rs}^{(d)}$ are given by

$$Q_{rs}^{(2)}=-\la_r^2+\la_r\la_s-\la_s^2,\quad
Q_{rs}^{(3)}=(\la_r^2+\la_s\la_k)(\la_s^2+\la_r\la_k), $$

$$Q_{rs}^{(4)}=-\g^{-2}(\la_r^2+\g^2)(\la_s^2+\g^2)
(\g^2+\la_r\la_k+\la_s\la_l)(\g^2+\la_r\la_l+\la_s\la_k),$$
where $\{ r,s,k,l\}=\{ 1,2,3,4\}$.

$$Q_{rs}^{(5)}=\g^{-8}(\g^2+\la_r\g+\la_r^2)(\g^2+\la_s\g+\la_s^2)
\ \prod_{k\neq r,s}(\g^2+\la_r\la_k)(\g^2+\la_s\la_k).$$
\v
$Proof.$ The statements are shown by straightforward computations.
We give some details for the interested reader.
For statement (a), observe that $P_1^{(d)}(B)$ is nonzero
only in the last row, where it coincides with the right-eigenvector
of $B$ for the eigenvalue $\la_1$,
 and that $P_d^{(d)}(A)$ is nonzero only in the last column,
where it coincides with the left eigenvector of $A$ for the eigenvalue
$\la_d$. Using the triangular form of the matrices, and the fact that
the nonzero diagonal entries are 
$P_1^{(d)}(\la_1)$ and  $P_d^{(d)}(\la_d)$, respectively, these
matrices can be computed easily.
It is also obvious that $P_1^{(d)}(B)P_d^{(d)}(A)$ is a multiple of
$E_{dd}$, and, multiplying it by $P_d^{(d)}(A)$ from the left,
that this multiple is equal to $Q_{1d}^{(d)}$.

Statement (b) follows from (a) and the fact that we obtain matrices in
ordered triangular form independent of the labeling, Prop. 2.4. Statement
(d) follows similarly from (a).
Statement (c) is straightforward.

\v\ni
{\bf 2.9. Main theorem} Let $K$ be an algebraically closed field.
(a)  Any simple $B_3$ module is uniquely determined by the eigenvalues of
$A$, up to a choice of a square root $\g^2$ of $\det(A)$ (for $d=4$) resp.
a 5-th root $\g$ of $\det(A)$ (for $d=5$).

(b) There exists a simple $B_3$ module $V$ of $K$-dimension
$d\leq 5$ if and only if the eigenvalues $\la_i$ of $A$
and the quantities $\g^2$ (for $d=4$) and $\g$ (for $d=5$), as defined in
(a), satisfy $Q_{rs}^{(d)}\neq 0$ for $r\neq s$ and $1\leq r,s \leq d$,
with  $Q_{rs}^{(d)}$ as defined in 2.8. The eigenvalues $\la_i$
need not be mutually distinct for this statement.

\v
$Proof.$ Let $V$ be a simple $B_3$-module.
By Propositions 2.5-7, we can assume a basis for $V$
such that $A$ and $B$ act via matrices as described there.
This shows part (a).

It remains to be shown for which values of the parameters the
representations given there are simple. Let us assume $Q_{ij}^{(d)}=0$
for some $i\neq j$, $1\leq i,j\leq d$. 
If $V$ is simple, we can assume
$P_i^{(d)}(A)\neq 0$, by the corollary of Prop. 2.2. 
Hence there exists
a vector $v$ for which $e_i=P_i^{(d)}(A)v$ is an eigenvector of $A$.
As $Q_{ji}^{(d)}=0$, also $P_j^{(d)}(B)P_i^{(d)}(A)=0$, by Prop.~2.8(b).
Obviously, $e_i$ is in the  subspace spanned by  
$\{ B^re_i, 0\leq r\leq d-2\}\cup \{ b_i\}$, which is $B$-invariant; it 
is a proper subspace of $V$ as the minimal polynomial of the restriction
of $B$ to it divides $P_j^{(d)}$.
We obtain the existence of a non-trivial $B_3$-submodule from Prop.~2.2,
i.e. $V$ can not be simple.

On the other hand, assume that $Q_{ij}^{(d)}\neq 0$ for all
$1\leq i\neq j\leq d$.
Let $W$ be a non-zero $B_3$-submodule of $V$. Then it contains at least one 
eigenvector of $A$, say $e_i$. As  $P_i^{(d)}(A)\neq 0$ by Prop. 2.8(a),
we can assume  this eigenvector to be of the form
$e_i=P_i^{(d)}(A)v$ for some $v\in V$, also if $A$ is not diagonalizable.
Let $j\neq i$. Then  $=ABAP_j^{(d)}(B)P_i^{(d)}(A)v$
is an eigenvector of $A$
with eigenvalue $\la_j$, provided it is nonzero. Using $(ABA)^{-1}P_i(B)ABA
=P_i(A)$, one easily checks that $P_i^{(d)}(B)e_j=Q_{ij}^{(d)}(ABA)e_i$,
which is nonzero by our assumptions. Hence $W$ contains eigenvectors
for each eigenvalue of $A$.

$Case$ 1: Assume $A$ has at least 2 distinct eigenvalues, which we label
$\la_1$ and $\la_2$. By triangularity of $A$, their eigenvectors
can be chosen in the form $(1,0,...)^t$ and $(*,1,0,...)^t$, which
are both in $W$. It is easy to see from this that
$e_1, e_2$, $b_1=ABAe_1$ and $b_2=ABAe_2$ together with $Ab_1$ (which
is the last column of $A$) span $V$, which therefore is equal to $W$.

$Case$ 2: Assume $A$ and $B$ has only one eigenvalue, say $\la$. Then
$P_1^{(d)}(A)P_d^{(d)}(B)=(A-\la)^{d-1}(B-\la)^{d-1}=Q_{1d}E_{dd}$.
As $Q_{1d}\neq 0$, the eigenspaces of $\la$ for both $A$ and $B$
are 1-dimensional. If $W$ is a non-zero $B_3$-submodule, it therefore
must contain the eigenvector $w_d$ of $B$. The set of vectors
$S=\{ (A-\la)^iw_d, i=0,1,\ ...\ d-1 \}$ is in $W$. As $(A-\la)^{d-1}\neq 0$,
so is its last column  $(A-\la)^{d-1}w_d$, by triangularity of $A$.
Hence $S$ is linearly independent. This finishes the proof.
\v\ni
{\bf Corollary} The simple $SL(2,\Z)$ modules and the simple  $PSL(2,\Z)$
modules are given by all simple $B_3$-modules for which
$\delta^2=1$ (for $SL(2,\Z)$) and for which $\delta =1$ (for  $PSL(2,\Z)$).
Observe that $\delta = -(\la_1\la_2)^3$ for $d=2$ and $\delta = 
 (\la_1\la_2\la_3)^2$ for $d=3$.
\v\ni 
{\bf 2.10 Parameter spaces for $B_3$ and $PSL(2, \Z)$} We define
ideals $I_d\subset K[\la_1, \ ...\ \la_d]$ (for $d\leq 3$)
and  $I_d\subset K[\la_1, \ ...\ \la_d, \g ]$ (for $d=4,5$) as follows:

\vl
$I_2=\langle \la_1^2-\la_1\la_2+\la_2^2\rangle$,
\vl
$I_3=\langle \la_i^2+\la_r\la_s,\  \{ i,r,s\} = \{ 1,2,3\}\rangle$, 
\vl
$I_4=\langle \{ \la_i^2+\gamma^2, i=1,2,3,4\}\  \cup 
\ \{ \gamma^2 +\la_i\la_j+\la_r\la_2,\  \{ i,j, r,s\} = \{ 1,2,3, 4\}\}
\ \cup\ \{ \g^4=\la_1\la_2\la_3\la_4\}\ \rangle$,
\vl
$I_5=\langle \{ \gamma^2+\gamma\la_i+\la_i^2,\  i=1,2,3,4\}\ \cup 
\ \{ \gamma^2+\la_i\la_j,\ 1\leq i<j\leq 5\}\ \cup\ \{ \g^5=
\la_1\la_2\la_3\la_4\la_5\} \rangle .$
\vl

Let $N_d$ be the zero set of $I_d$ in $K^d$ (resp. in $K^{d+1}$
for $d>3$). Observe that
$N_d$ is invariant under the action of the symmetric group,
acting via permuting the coordinates $\la_i,\ 1\leq i\leq 5$
in $K^d$. The main theorem can now be reformulated as follows:
\v\ni
{\bf Main Theorem$'$}
There exists a 1-1 correspondence
between equivalence classes of simple $B_3$ modules of dimension $d\leq 5$
and the $S_d$ orbits in $K^d\backslash N_d$.
\v\ni
{\bf 2.11 Remarks} 1. After this research was completed we became aware of
a number of related approaches and results. We learned
from P. Deligne that results similar to ours are obtained in
N. Katz's work [K] on rigid local systems by a more general and less
elementary method.
In particular, he sketched to us how one can obtain the classification
of 5-dimensional braid representations from Katz's work (see also Remark 4
below).

2. We also learned of another approach using  quiver theory (see [Ws]
and, for a more general  approach, [S]). The following result
appears explicitly in a preprint by Westbury:

Recall that in the $PSL(2,\Z)$-quotient of $B_3$ the elements 
$ABA$ and $AB$ have order 2 and 3 respectively.
Let $n_1$ and $n_2$ be the dimensions of the eigenspaces of $ABA$,
and let $m_1$, $m_2$ and $m_3$ be the dimensions of the eigenspaces
of $AB$. Then there exists an indecomposable representation of  $PSL(2,\Z)$
if and only if
$n_i\geq m_j$ for $i=1,2$ and $ j=1,2,3$.
Moreover, the parameter space for representations of this type has 
dimension
$d^2 - n_1^2 - n_2^2 -m_1^2 - m_2^2 - m_3^2 + 1$.
This result is in accordance with our findings for $d\leq 5$, but does 
not say anything about when the representations are simple (e.g.
for $d=4$ there also exist indecomposable but not simple representations
with $AB$ only having 2 eigenvalues).

3. We do not expect our methods to work for dimension $\geq 6$ without
significant changes. For one, we can not expect ordered triangular forms
for $A$ and $B$: Indeed, if this were the case, $Tr(ABA)=0$
for any simple $6$-dimensional representation of $B_3$, by Prop. 1.3(a),
i.e. $n_1=n_2=3$. To get a counterexample, it suffices to find
$6\times 6$ matrices of order 2 and 3 which generate the whole ring of
$6\times 6$ matrices and such that $n_1=4$ and $n_2=2$.

Similarly, the result quoted under 2 also shows that we will have
more parameters than the eigenvalues for $d\geq 6$ and suitable
values of $m_1, m_2, m_3, n_1, n_2$.

4. There exists a beautiful and
simple argument which considerably narrows down which polynomials can occur in
the simplicity statement of Theorem 2.9, and which also explains
the nature of the factors in the polynomials $Q_{ij}^{(d)}$ to some extent.
We were told this argument by P. Deligne:

Let $V$ be a $d$-dimensional $B_3$ module on which the central element
$\zeta=(\sigma_1\sigma_2)^3$ acts
as a scalar $\delta$.
Assume that $V$ has an $r$-dimensional $B_3$-submodule $W$
on which $\sigma_1$ acts with eigenvalues $\la_1,\ ...\ \la_r$.
Comparing determinants of the restriction of $\zeta$ to $W$, we get that 
$$(\la_1\ ...\ \la_d)^{6r}=(\la_1\ ...\ \la_r)^{6d}.\eqno(*)$$
Hence whenever $\zeta$ acts as a scalar and $(*)$ is NOT satisfied,
the representation has to be simple.
Unfortunately, this argument produces a sufficient but not a necessary
condition for $V$ being simple.

5. A result which could be related to our findings is
Coxeter's classification of  finite quotients of  $B_3$ defined by
the additional relation $\sigma_i^p=1$, for $i=1,2$: the additional
relation defines a finite quotient if and only if $p\leq 5$ (see [C]).
\vskip .6cm
\ni
{\bf 3. An Application to Tensor Categories}
\v\ni
{\bf 3.1} We shall use the results of the previous section to 
compute categorical dimensions for certain simple objects in
braided tensor categories related to exceptional 
Lie groups. We will not give much background
information  about tensor categories here (see  e.g. [DM], [JS],[KW] or [T]),
as the application itself is rather elementary.
The reader not  familiar
with the categorical language should think of the categories as
representation categories of quantum groups.

Let $\Ca$ be a semisimple ribbon tensor category with trivial object
$\1$ over an algebraically closed field $K$ (see [T] for
precise definitions);  this
means, in particular,  that $End(X)$ is a semisimple $K$-algebra 
for any given   object $X$ in $\Ca$, and the homomorphisms between 2
objects in $\Ca$ form a vector space over $K$.

In the following, let $Z$ be a simple selfdual object in $\Ca$;
this means that the trivial object $\1$ appears with multiplicity
one as a direct summand in $Z^{\otimes 2}$. 
Let $p\in \End(Z^{\otimes 2})$
be the projection onto $\1$, and let $p_1=p\otimes 1$ and $p_2=
1\otimes p$ be elements in $\End(Z^{\otimes 3})$. 
In the following we will always assume the following (rigidity) condition
$$p_2p_1p_2\neq 0\eqno(3.1)$$
As $Z\otimes \1$
is canonically isomorphic to the simple object $Z$, 
$p_2(a\otimes 1)p_2$ is equal
to a scalar multiple of $p_2$ for any $a\in \End(Z^{\otimes 2})$.
Let $Y$ be a direct summand in $Z^{\otimes 2}$, with respect to
a chosen direct sum decomposition
of $Z^{\otimes 2}$. This decomposition
defines a projection $p_Y$ onto $Y$. We define the
categorical dimension $\dim Y$ by
$$(\dim Y)p_2=(\dim Z)^2p_2(p_Y\otimes 1)p_2.\eqno(3.2)$$
 This definition depends on a choice
of $\dim Z$. In the following we will always require  $\dim\ \1 = 1$.
Taking for $p_Y$ the projection $p_1$, it follows from (3.1) and (3.2)
that $\dim Z$ is uniquely determined up to a sign;
$\dim Y$ is independent of the choice of the sign.
So, in particular, if $Z^{\otimes 2}$ contains a subobject 
isomorphic to $Z$, $\dim Z$ is uniquely determined by (3.2).
It is not hard to check that for subobjects of
$Z^{\otimes 2}$ this definition  coincides with the
usual definition of categorical trace as e.g. for ribbon tensor
categories (see [T]).
\v\ni
{\bf 3.2} Assume that $\Ca$ allows a braiding (see e.g. [JS], [T]). 
For the purpose of this article,
it suffices to know that this entails that there exists a canonical
endomorphism $c\in \End(Z^{\otimes 2})$ such that $c_1c_2c_1=c_2c_1c_2$,
with $c_i$ elements in $\End(Z^{\otimes 3})$ defined as before the $p_i$.
We have the following simple 
\v\ni
{\bf Lemma} Assume that $Z^{\otimes 2}$ decomposes as a direct
sum $\bigoplus_i Y_i$ of $d$ mutually nonisomorphic simple objects $Y_i$, 
each of which has non-zero dimension.
Moreover, we assume that $c$ acts on $Y_i$ via the scalar $\la_i$,
and $\la_i\neq \la_j$ for $i\neq j$.

Then $\Hom(Z,Z^{\otimes 3})$ is a simple $d$-dimensional  $B_3$-module,
with the action defined via $\sigma_jf=c_j\circ f$ for all
$f\in \Hom(Z,Z^{\otimes 3})$ and $j=1,2$; moreover, each eigenvalue of
$c_j$ has multiplicity 1.
\v
$Proof$. Let $p^{(i)}\in \End(Z^{\otimes 2})$ be the projection onto $Y_i$,
 with $p=p^{(1)}$ being
the projection onto $\1\subset Z^{\otimes 2}$ and let $\iota: \1 \to 
Z^{\otimes 2}$ be a nonzero homomorphism. 
It follows from rigidity that $(p^{(i)}\otimes 1)\circ (1\otimes \iota)$
is a nonzero homomorphism and that $\dim \Hom(Z, Z^{\otimes 3})=
\dim \End(Z^{\otimes 2})=d$; it can also be checked explicitly for the examples
below. Hence $V=\Hom(Z,Z^{\otimes 3})$ is a $d$-dimensional
vector space on which both $c_1$ and $c_2$ act via concatenation of
morphisms. Moreover, it has a basis of eigenvectors 
$(p^{(i)}\otimes 1)\circ (1\otimes \iota)$ of $c_1$.
If this representation of $B_3$ were not simple, fix
a composition series of $V$, and pick an eigenprojection $p^{(i)}$ of $c$
such that the eigenprojections $p^{(i)}\otimes 1$ and $p_1=p\otimes 1$ of $c_1$
act nonzero on different simple factors. As $p_1$ is conjugate
to $p_2=1\otimes p$ (as $c_1$ and $c_2$ are conjugate), we can conclude
$p_2(p^{(i)}\otimes 1)p_2=0$, 
which would contradict the assumption
about nonzero dimensions.
\v\ni
{\bf Corollary} Using the notations of the Lemma, with $Y_1\cong \1$, we have
$$\dim Y_i={Q_{1i}^{(d)}(\dim Z)^2\over P_1^{(d)}(\la_1)P_i^{(d)}(\la_i)},$$
where $Q_{1i}$, $P_1^{(d)}(\la_1)$ and $P_i^{(d)}(\la_i)$ are as in 
Section 2.8.
\v\ni
$Proof.$ Observe that the projection onto the eigenspace of $A$
with eigenvalue $\la$ is given by $P_i^{(d)}(A)/P_i^{(d)}(\la_i)$
The claim now follows from (3.2) and Prop. 2.8(a).
\v\ni
{\bf 3.3} 
Let us first consider the case with $\Ca$ a
braided tensor category whose Grothendieck semiring is isomorphic
to the semiring of the representation category of the Lie group $G$ with
$G$ being an orthogonal group $O(N)$ or a
symplectic group $Sp(N)$. We take as object $Z$ the vector
representation of $G$. It is well-known that in this case
$Z\otimes Z\cong \1\oplus X\oplus Y$, with $X$ and $Y$ simple objects
corresponding to the antisymmetrization and a subrepresentation
of the symmetrization
(traceless tensors) of  the 2nd tensor power of
the vector  representation.
It can be shown
that in this case (only using the assumptions of braiding, and 
the given Grothendieck semi-ring), the eigenvalues are of
the form $\alpha q, -\alpha q^{-1}, \alpha r^{-1}$, with $\alpha$
a 4-th root of unity and $q,r\in K$.
Moreover, $r$ is $\pm$ a power of $q$, depending
on the given category (see below). Using 
the notations $[n]=q^n-q^{-n}$ and $[\la +n]=rq^n-r^{-1}q^{-n}$,
it follows from Corollary 3.2 that
$$\dim Z=\alpha^2({[\la]\over [1]}+1)=\alpha^2({r-r^{-1}\over q-q^{-1}}+1),$$
and
$$\dim X={[\la -1]+[2]\over [2]}\ {[\la]\over [1]},
\quad
\dim Y={[\la +1]+[2]\over [2]}\ {[\la ]\over [1]}.$$
This method can be extended to define similar functions also
for objects in higher tensor powers of $Z$, which are labeled by
Young diagrams (see [Wn, Theorem 5.5]). The $q$-dimensions
for $G=O(N)$ are obtained by setting  $\la = N-1$, for $G=Sp(N)$
they are obtained by setting  $\la=-N-1$.
More generally, it is possible to reconstruct such categories similarly
as it was done for categories of type A in [KW]; however, in this case
we would need the assumption of the category being braided.
So in this comparatively simple case knowledge about
braid representations allows us to reconstruct the tensor category.
\v\ni
{\bf 3.4} The discussion in 3.3 can be carried over to exceptional
Lie groups to some (so far rather limited)  extent.
We use the notations $[n]=s^{n/2}-s^{-n/2}$
and $[\la+n]=t^{1/2}s^{n/2}-t^{-1/2}s^{-n/2}$ for $n\in \N$ and $\la$ a formal 
variable; in our formulas, the result will be independent of the choice
of square roots of $s$ and $t$. 
It can be shown that for an exceptional Lie algebra $\go$,
the 2nd tensor product $\go^{\otimes 2}$ of its adjoint representation
decomposes as a direct sum $\1\oplus\go\oplus X_2\oplus Y_2\oplus Y_2^*$.
Here we use the same notations as in [De]. 
We are going to use the following results about quantum groups,
which are essentially due to Drinfeld [Dr2] (for a description of
the relationship between operators $c$ (braiding operators) and 
quantum Casimir (= twist) also see e.g. [T]); here
$a_Y$ denotes the scalar via which the Casimir acts on a simple 
$\go$-module $Y$:
\vl
(i) the braiding operator $c$ acts on the simple summand 
$Y\subset \go^{\otimes 2}$ via the scalar $\pm q^{a_Y-2a_\go}$, with the sign
depending on whether $Y$ is in the symmetrization or antisymmetrization
of $\go^{\otimes 2}$,
\vl
(ii) the operator $(c_1c_2)^3$ acts on the simple summand 
$K\subset \go^{\otimes 3}$ via the scalar $q^{2a_K-6a_\go}$.
\vl
Now setting $q^{2a_\go}=s^{-3}$ and $q^{a_{Y_2}-2a_\go}=t$, and using
the formulas for the action of the Casimir in [De], one sees easily that
the eigenvalues  of $c$ (acting on $\go^{\otimes 2}$)
are $s^6$ (for $\1$), $-s^{3}$ (for $\go$), $-1$ (for $X_2$),
$t$ (for $Y_2$) and $st^{-1}$ (for $Y_2^*$). It  follows from this
that  det$_V(c_1)=s^{10}$, where  det$_V(c_1)$ is the determinant of the linear
operator via which $c_1$ acts on $V$.

In order to completely determine the 5-dimensional  representation of $B_3$
on $V=\Hom(\go, \go^{\otimes 3})$, we only need to compute the scalar $\delta$
by which the central element $(c_1c_2)^3$ acts on $V$. This is equivalent to
determining the 5-th root $\g$ of det$_V(c_1)$ given by
$\delta = \g^6=\det_V(c_1)\g = s^{10}\g$. By (ii),
we get $\delta = q^{2a_\go - 6a_\go}=s^{12}$. Combining the last 2 formulas,
we get $\g=s^2$. One can now check, using Corollary 3.2, that
$$\dim \go =
{[4][\la -6][\la +5]\over [2][\la -1][\la]},$$

$$\dim {X_2}=
{[5][\la -6][\la +5][\la -4][\la +3] [2\la +4][2\la -6]
\over
[1][\la ][\la -1][\la +2][\la - 3][2\la][2\la -2]},$$

$$\dim {Y_2}=
{[6][5][4][\la +5][\la-4][3\la -6] \over
[2][\la -1][\la][2\la][2\la -1][\la -2]},$$

$$\dim {Y_2^*}=
{ [6]][5][4][\la -6][2\la+6][3\la +3] [\la +2]\over
[2][\la -1][\la][2\la-2][2\la -1][2\la +4][\la +1]}.
$$
These rational functions are $s$-deformations
of  formulas given in [De] for the classical  case.
\v\ni
{\bf 3.5} The computations above in the orthogonal case
were used by Toledano Laredo in [TL] in the
context of fusion of representations of loop groups,
based on unpublished notes of the second named author.
We have been informed
by Antony Wassermann that he uses our results for 5-dimensional 
braid representations in connection with fusion of representations of
loop groups corresponding to exceptional Lie groups.

\vskip 1.5cm
\ni
{\bf References}
\v
\item{[BLM]} J.S Birman, D. D. Long and J. A. Moody
Finite-dimensional representations of Artin's braid group,
Contemp. Math., 169,
Amer. Math. Soc., Providence, RI, 1994, 123--132.
\item{[C]} H.S.M. Coxeter, Factor groups of the braid group,
Proceedings of the Fourth Can. Math. Congress, Banff 1957,
University of Toronto Press (1959), 95-122
\item{[De]}  P. Deligne, La s\' erie exceptionnelle de groupes de Lie. 
C. R. Acad. Sci. Paris Sér. I Math. 322 (1996), no. 4, 321--326. 
\item{[DM]} P. Deligne and J. Milne, Tannakian categories, Lecture
Notes in Mathematics, vol 900, Springer Verlag
\item{[Dr1]} V. Drinfeld, Quantum groups, Proceedings for the ICM Berkeley
1986 798-820.
\item{[Dr2]} V. Drinfeld, On almost cocommutative Hopf algebras,
 Leningrad Math. J. 1 (1990) 321-243.
\item{[JS]} A. Joyal and R. Street, The geometry of tensor calculus,
Adv. Math. 88 (1991)
\item{[K]} N. Katz, Rigid local systems, Annals of Mathematics Studies,
139, Princeton University Press
\item{[KW]} D. Kazhdan and H. Wenzl, Reconstructing monoidal categories,
Advances in Soviet Math. Vol 16 Part 2 (1993) 111-136
\item{[S]} A. Schofield, General representations of quivers. 
Proc. London Math. Soc. (3) 65 (1992), no. 1, 46--64. 
\item{[TL]} V. Toledano Laredo, Fusion of positive energy representations
of LSpin$_{2n}$, Thesis, University of Cambridge, 1997
\item{[T]} V.G. Turaev, Quantum Invariants of Knots and
$3$-Manifolds, de Gruyter (1994)
\item{[V]} P. Vogel, Algebraic structures on modules of diagrams, preprint
\item{[Wa]} A. Wassermann, Operator algebras and conformal field
theory III, Invent. Math. (1998)
\item{[Wn]} H. Wenzl, Quantum groups and subfactors of type B, C and D, 
Comm. Math. Phys. 133 (1990) 383-432
\item{[Ws]} B. Westbury,  On the character varieties of the modular group,
preprint, University of Nottingham, 1995
\v\v
\ni Department of Mathematics

\ni University of California at San Diego

\ni La Jolla, CA 92093-0112, USA
\v\ni
email: {\tt ituba$\char'100$euclid.ucsd.edu} and 
{\tt wenzl$\char'100$brauer.ucsd.edu}

\end

\item{[Ji]} M. Jimbo, Quantum $R$-matrix for the generalized Toda
system, Comm. Math. Phys. 102,4 (1986), 537-547.

\item{[Kc]} V. Kac, Infinite dimensional Lie algebras, Cambridge University
Press

If a matrix $A$ has eigenvalues $\la_i$, $1\leq i\leq d$, each with
multiplicity 1, then the projection $E_r(A)$ onto the eigenspace of $\la_r$,
with kernel being the span of the other eigenspaces,
is given by the formula
$E_r(A)=P_r^{(d)}(A)/P_r^{(d)}(\la_r)$,
with $P_r^{(d)}$ as in 2.8. If $A$ and $B$ are matrices 
in a simple $d$-dimensional representation of $B_3$, we get
$$E_s(B)E_r(A)E_s(B)= {Q_{rs}\over P_s^{(d)}(\la_s)P_r^{(d)}(\la_r)}.$$

For the group $PSL(2,\Z)$, the zero set $N_d$ can be described more
explicitly, stemming from the condition that $(AB)^3=\det(A)^{6/d}1$ has to be
equal to 1. In the following we always assume $\theta$ to be a primitve
6-th root of unity, i.e. $\theta = e^{\pm \pi i/3}$.

$d=2$:  One checks by direct computation with the matrices
that $(AB)^3=-(\la_1\la_2)^31$; hence we get $\la_1\la_2=-\theta^{2j}=
-\theta^{-j}$, for $j=0,1,$ or 2.
Observe that the condition for $N_2$ can be rephrased as $\la_1/\la_2\neq
\theta^{\pm 1}$. This can be combined to saying that there exists  a
simple 2-dimensional $PSL(2,\Z)$-representation if and only if
$$\la_1\la_2=-\theta^{2j}, \quad j=0,1,2\quad {\rm and}
\quad \la_1\neq \pm\theta^j.$$

$d=3$: Here one checks similarly as for $d=2$ that there exists a simple
3-dimensional  $PSL(2,\Z)$-representation if and only if
$$\la_1\la_2\la_3=\pm 1,\quad \la_i^3\neq  -(\pm 1).$$

$d=4$: Let $\g=\det(A)^{1/2}$. Then $\g^3=1$ from which we get $\g=
\theta^{2j}$, with $j=0,1,2$. The equation
$ \gamma^2 +\la_i\la_j+\la_r\la_2\neq 0$ is equivalent to
$\la_i\la_j/\g^2\neq \theta^{\pm 2j}$.
Hence  there exists a simple
4-dimensional  $PSL(2,\Z)$-representation if and only if
$$\la_i\la_j\neq \theta^{2i+2j},\ i\neq j\quad {\rm and}
\quad \la_i^2\neq-\theta^2.$$

$d=5$: Let $\g=\det(A)^{1/5}$. Then $\g$ has to be a
6$^{th}$ root of unity, not necessarily primitive, and we get the
following necessary and sufficient conditions for the existence
of a simple 5-dimensional  $PSL(2,\Z)$-representation:
$$\la_i\neq \g e^{\pm 2\pi i/3}\quad {\rm and}\quad
\la_i\la_j\neq -\g^2\ {\rm for}\ i\neq j.$$
\v\ni

2. One of the motivations of our study was to use our results in order
to study braided tensor categories. In [KW], a detailed knowledge
of braid representations into Hecke algebras enabled the authors to
completely classify tensor categories whose Grothendieck ring coincides 
with the
one of the representation category of $SU(k)$ or one of 
the associated fusion categories. It is
fairly straightforward to  extend this approach to braided categories
whose Grothendieck ring coincides with the one of the representation
category of the orthogonal group $O(n)$. This was done in unpublished
notes by the second named author, and some of the computations
occurring there  were subsequently used in [TL].
This only involved representations of $B_3$ of dimension $\leq 3$.
To study similar questions in connection with exceptional Lie algebras,
representations of $B_3$ will appear up to dimension 5, but not more
(see [De]).
We plan to publish applications of our results for that somewhere else.
We would also like to mention that
Antony Wassermann has informed us that he can give an operator
algebraic proof of the Verlinde formula for $E_8$ using our results.

$$A=
\left(
\matrix{\la_1&*&&...&&*\cr 0&\la_2&&...&&*\cr&&&&&\cr 0&&...&&0&\la_d\cr}
\right),\quad
B=
\left(
\matrix{\la_d&0&&...&&0\cr *&\la_{d-1}&&...&&0\cr&&&&&\cr *&&...&&*&\la_1\cr}
\right).$$

(a)\quad   $ABA$ is a skew diagonal matrix, with $(ABA)_{i,d+1-i}=
(-1)^{i+1}\delta^{1/2}$,
\vskip .1cm
(b)\quad  $b_{ij}=(-1)^{i+j}a_{d+1-i,d+1-j}$,
\vskip .1cm
(c)\quad  $\sum_{k=1}^i(-1)^{k+1}a_{d+1-i,d+1-k}a_{k,d+1-i}=
\delta^{1/2}/\la_i$,
\vskip .1cm
(d)\quad  $(BA)_{i,d+1-i}=(-1)^{i+1}\delta^{1/2}/\la_i$, 
\vskip .1cm

(e)\quad  If $i+j>d+1$, $\sum_{k=1}^{min(i,j)}(-1)^{k+i}a_{d+1-i,d+1-k}a_{k,j}
=0$

$claim$: If $w_2=e_2-(c_{d2}/c_{d1})e_1\in W=span\{ e_1,e_2\}\cap
span\{ b_1,b_2\}$, then $V$ has a nontrivial $B_3$-submodule.

If dim$W=2$, $W=span\{ e_1,e_2\}=span\{ b_1,b_2\}$ 
would be a subspace invariant under both $A$
and $B$. Hence we can assume  dim$W=1$, and therefore $W=span\{ w_2\}$.
As $ABAW=W$, we get $ABAw_2=\mu w_2$ for some scalar $\mu$, and therefore
$$BAw_2=\mu A^{-1}w_2.$$
Using the definition of $w_2$, and the fact that the $e_i$'s are eigenvectors
of $A$, the last equation can be transformed into
$$B((\la_2-\la_1)e_2+\la_1w_2)=\mu(\la_2^{-1}-\la_1^{-1})e_2
+\mu\la_1^{-1}w_2.$$
As $w_2\in$ span$\{ b_1,b_2\}$, we get $Be_2\equiv \gamma e_2$
mod span$\{ b_1,b_2\}$, for some scalar $\gamma$, i.e. $e_2\equiv b_i$
mod  span$\{ b_1,b_2\}$ for some $i>2$.
If $i\neq d$, then $c_{d2}=0$(!), and therefore $e_2=w_2=\mu^{-1}ABAw_2=
\mu^{-1}ABAe_2=\mu^{-1}b_2$.
Hence span$\{ e_2\}$ is an invariant subspace.

So we can assume $e_2\equiv b_d$ mod span$\{ b_1,b_2\}$, and $c_{d2}\neq 0$.
This implies 
$$e_1={c_{d1}\over c_{d2}}(e_2-w_2)\in span\{ b_1,b_2,b_d \}.$$
As $w_2\in span\{ b_1,b_2 \}$, we also get $e_2\in span\{ b_1,b_2,b_d \}$.
Solving for $b_d$ in the expansion of $e_1$ in terms of $b_i$'s
(possible, as $c_{d1}\neq 0$), we get from the above that $b_d\in
span\{ b_1,b_2,e_1\}$. Multiplying this with $ABA$, it follows that
$e_d\in span\{ e_1,e_2, b_1\}\subset  span\{ b_1,b_2,b_d \}$.
Hence $ span\{ e_1,e_2, e_d\}=span\{ b_1,b_2,b_d \}$ is a $B_3$-submodule,
which finishes the proof of our claim.

Then claim (a) will follow
as soon as we have shown that
$\{ w,\ ABw,\ (AB)^2w\}$ spans $V'$. Indeed, as $(AB)^3w$ is a multiple of
$w$, it follows that $V'$ is invariant under both $AB$ and $ABA$.

Assume there exist vectors $w_i,\ i=1,2,\ ...\ d$ satisfying
\vskip .1cm
(i) $w_i\in span \{ e_1,\ ...\ e_i\}\cap  span \{ b_{d+1-i},\ ...\ b_d\}$
\vskip .1cm
(ii) If $w_i=\sum_{j=1}^i \alpha_{ij}e_j=
\sum_{j=d+i-i}^d \beta_{ij}b_j$, then $\alpha_{ii}
\neq 0\neq \beta_{\bar i,\bar i}$.
\vskip .1cm
\ni Then  $\{ w_1,\ ...\ w_d\}$  is a basis with respect to which
 $A$ and $B$ have ordered triangular form.

\v
$Proof.$ By condition (i) and (ii), $\{ w_1,\ ...\ w_d\}$ is a basis.
We define subspaces
$$W_i= \span \{ w_1,\ ...\ w_i\} =  \span \{ e_1,\ ...\ e_i\}\eqno(2.1.1)$$ and
$$ W_i'=  \span \{ w_{d+1-i},\ ...\ w_d\} =  
\span \{ b_{d+1-i},\ ...\ b_d\}.\eqno(2.1.2)$$
Then $AW_i\subset W_i$ and $BW_i'\subset W_i'$, i.e. $A$ and $B$ 
have triangular shapes as required. The diagonal entries can easily
be computed by induction on $i$ by taking into account $\det(A_{|W_i})=
\prod_{j=1}^i \la_j$ and $\det(B_{|W'_i})=
\prod_{j=d+1-i}^d \la_j$.
\v\ni
{\bf 2.2. Lemma} For any 2 generalized eigenvectors $e_1$ and $e_2$,
we can find a generalized eigenvector, labeled $e_d$ such that

(a) $c_{d1}\neq 0$, and the set $\{w_1, w_d\}$ is linearly independent
for $w_1=e_1, w_d=b_1$,

(b) $w_2=e_2-(c_{d2}/c_{d1})e_1\in span\{ e_1,e_2\}\cap  
span\{ b_1,\ ...\ b_{d-1}\}$,

(c) $w_{d-1}=ABAw_2\in  span\{ e_1,\ ...\ b_{d-1}\}\cap  
span\{ b_1, b_2\}$.

(d) If $W=\span\{e_1,e_2\}\cap\span\{w_1,w_2\}\neq 0$, then $V'=
\span\{e_1,e_2,w_1,w_2\}$ is a 3-dimensional $B_3$-invariant subspace. 
\v
$Proof.$ If $e_1\in$ span$\{b_1\}$, 
span$\{ e_1\}$= span$\{ b_1\}$ would be a subspace
both invariant under $A$ and $B$, contradicting simplicity of $V$.
Hence there exists a $j\neq 1$ such that $c_{j1}\neq 0$, which we  relabel
to be $d$. Statements (b) and (c) are easy consequences of this and
Lemma 2.1. 

To prove (d), it suffices to show that $V'$ is invariant
under both $ABA$ and $AB$. Invariance under $ABA$ follows from
Lemma 1.2(d) and (b). If $\dim W = 2$, then $W=
\span\{e_1,e_2\}=\span\{w_1,w_2\}$ obviously is invariant under both
$A$ and $B$.  If $\dim W = 1$, let $w=\alpha e_1+\beta e_2$ be a
basis vector for it. Again using Lemma 1.2, we see that $W$ is invariant
under $ABA$, and $ABAw=\alpha b_1+\beta b_2=\mu w$ for some scalar $\mu$.
If $e_1$ and $e_2$ are eigenvectors with distinct eigenvalues $\la_1$
and $\la_2$, we  compute
$$ABw=B^{-1}BABw=\la_2^{-1}w+\alpha(\la_1^{-1}-\la_2^{-1})b_1\quad
{\rm and}\quad (AB)^2w= \mu\la_2w+\mu\alpha(\la_1-\la_2)e_1.$$
We can assume $\alpha\neq 0$ (otherwise $W$ would be the eigenspace for
$\la_2$ for both $A$ and $B$). It is now easy to check that
both $\{ e_1, b_1, w\}$ $\{ w, ABw, (AB)^2w\}$ span $V'$.
Similarly, if $e_1$ is an eigenvector with eigenvalue $\la$ and
$(A-\la)e_2=e_1$, one computes
$$ABw=\la^{-1}w-\la^{-2}e_1\quad {\rm and}\quad
(AB)^2w=\mu w+\mu\beta e_1.$$
>From this one deduces the claim as before.

{\bf 2.4. Lemma} There are bases  $\{ w_1,\ ...\ w_d\}$
for $d=4,5$ for which $A$ and $B$ are in ordered triangular form.
\v
$Proof.$ Construct $w_1,w_2, w_d$ as in Lemma 2.2. We also know from this
lemma that $e_1, e_2, b_1, b_2$ are linearly independent.
We can now easily check that the conditions of Lemma 2.1 are satisfied
for $d=4$, with $w_3=ABAw_2$. This proves the lemma for $d=4$,

To prove the lemma for $d=5$, we can again assume $w_2\not\in
span\{ b_1,b_2\}$, by Lemma 2.2. Hence we can assume
$w_2=\sum \beta_{2j}b_j$ with $\beta_{24}\neq 0$, after possible
relabeling. By our assumptions in $w_1$ and $w_2$ we can find
$\delta_1$ and $\delta_2$ such that
$$w_3=e_3-\delta_2w_2-\delta_1w_1 \subset span\{b_1,b_2,b_3\}.$$
Then $w_i=e_i+\sum_{j<i} \alpha_{ij}e_j$ and hence $\{ w_i\}$ is
linearly independent. The conditions of Lemma 2.1 can now again
be checked easily.
\v\ni